\documentclass[letterpaper]{article}
\usepackage{proceed2e}
\usepackage[margin=1in]{geometry}

\usepackage{times}

\usepackage{amsbsy}
\usepackage{amsmath}
\usepackage{amsfonts}
\usepackage{amsthm}
\usepackage{amssymb}
\usepackage{algorithmic}
\usepackage{algorithm}
\usepackage{enumerate}
\usepackage{graphicx}
\usepackage{multirow}
\usepackage{natbib}
\setcitestyle{round}
\usepackage{subfigure}
\usepackage{listings}
\usepackage{xcolor}
\theoremstyle{plain}
\newtheorem{theorem}{Theorem}

\newtheorem{lemma}{Lemma}
\newtheorem{condition}{Condition}
\theoremstyle{definition}
\newtheorem{experiment}{Model} 
\theoremstyle{remark}
\newtheorem*{remark}{Remark}

\usepackage[normalem]{ulem}

\title{Testing for Conditional Mean Independence with Covariates through Martingale Difference Divergence}

\author{ {\bf Ze Jin\thanks{\, Corresponding author. Email address: zj58@cornell.edu.}} \\
Department of Statistical Science \\
Cornell University\\
Ithaca, NY 14850 \\
\And
{\bf Xiaohan Yan}  \\
Department of Statistical Science\\
Cornell University          \\
Ithaca, NY 14850 \\
\And
{\bf David S. Matteson\thanks{\, Research support from an NSF Award (DMS-1455172), a Xerox PARC Faculty Research Award,
  and Cornell University Atkinson Center for a Sustainable Future (AVF-2017).}}   \\
Department of Statistical Science\\
Cornell University \\
Ithaca, NY 14850    \\
\\
}

\begin{document}

\maketitle

\begin{abstract}
As a crucial problem in statistics is to decide whether additional variables are needed in a regression model.
We propose a new multivariate test to investigate the conditional mean independence of $Y$ given $X$ conditioning on some known effect $Z$,
i.e., $\textrm{E}(Y | X, Z) = \textrm{E}(Y | Z)$.
Assuming that $\textrm{E}(Y | Z)$ and $Z$ are linearly related, we reformulate an equivalent notion of conditional mean independence through transformation,
which is approximated in practice. We apply the martingale difference divergence \citep{shao2014martingale} to measure conditional mean dependence,
and show that the estimation error from approximation is negligible, as it has no impact on the asymptotic distribution of the test statistic under some regularity assumptions.
The implementation of our test is demonstrated by both simulations and a financial data example.
\end{abstract}


\section{INTRODUCTION}\label{intro}


Testing (conditional) dependence and conditional mean dependence
plays an important role in statistics with various applications,
including variable selection \citep{szekely2014partial, park2015partial}
and feature screening \citep{li2012feature, shao2014martingale}.
Both areas attracted tremendous attention in the last two decades,
as datasets have increased in size and dimension.
Let $X \in \mathbb{R}^p$, $Y \in \mathbb{R}^q$, $Z \in \mathbb{R}^r$ be the three random vectors of interest,
and denote pairwise independence by ${\perp\!\!\!\!\perp}$.

Measures of (conditional) dependence have been extensively studied.
\cite{szekely2007measuring} proposed distance covariance (dCov) to capture the non-linear and non-monotone
pairwise dependence between $X$ and $Y$,
and dCov = 0 if and only if pairwise independence ($X {\perp\!\!\!\!\perp} Y$) holds.
\cite{jin2017generalizing} extended distance covariance to mutual dependence measures (MDMs).
To capture the conditional dependence between $X$ and $Y$ given $Z$,
\cite{szekely2014partial} generalized distance covariance to partial distance covariance (pdCov),
however, pdCov = 0 is not equivalent to conditional independence ($X {\perp\!\!\!\!\perp} Y | Z$).
\cite{wang2015conditional} extended distance covariance to conditional distance covariance (CDCov)
using kernel estimators, and CDCov = 0 if and only if conditional independence holds.
Under a linear assumption between $X, Y$ and $Z$,
\cite{fan2015conditional} converted testing conditional independence to testing independence,
and applied distance covariance to measure the dependence of estimated variables.
Moreover, conditional dependence is known as Granger causality in time series analysis.
\citet{hiemstra1994testing}, \citet{su2007consistent}, and \citet{chen2012testing}
each introduced non-parametric tests for non-linear Granger causality based on conditional probabilities and characteristic functions.

Likewise, various measures of conditional mean dependence
have been broadly developed as well.
Testing the conditional mean independence of $Y$ given $X$, i.e.,
\begin{equation}\label{h1}
H_0: \textrm{E}(Y | X) = \textrm{E}(Y) \,\,\,\, a.s., \quad H_A: o.w.
\end{equation}
provides insight on whether $X$ contributes to the conditional mean of $Y$.
\cite{shao2014martingale} generalized distance covariance to martingale difference divergence (MDD),
and MDD = 0 if and only if (\ref{h1}) holds.
Testing the conditional mean independence of $Y$ given $X$ conditioning on some known effect $Z$, i.e.,
\begin{equation}\label{h2}
H_0: \textrm{E}(Y | X, Z) = \textrm{E}(Y | Z) \,\,\,\, a.s., \quad H_A: o.w.
\end{equation}
sheds light on whether $X$ contributes to the conditional mean of $Y$ 
when taking known dependence on $Z$ into account.
\cite{park2015partial} generalized martingale difference divergence to partial martingale difference divergence (pMDD),
however, pMDD = 0 is not equivalent to (\ref{h2}).
\cite{fan1996consistent}, \cite{lavergne2000nonparametric}, and \cite{ait2001goodness}
each introduced non-parametric tests for (\ref{h2}) using kernel estimators of conditional expectations.
Assuming a linear model between $Y$ and $(X, Z)$,
\citet{lan2014testing} generalized the classical partial F-test \citep{chatterjee2015regression}
to a partial covariance-based (pcov) test for (\ref{h2}) in the high-dimensional setting,
and \citet{tang2017testing} further proposed a hybrid test for (\ref{h2})
through finding the most predictive covariate based on both maximum-type and sum-type statistics.
Conditional mean independence conditioning on covariates is known as Granger causality \textit{in mean} in time series analysis.
\citet{raissi2011testing} proposed a parametric test for linear Granger causality in mean based on vector autoregressive (VAR) models,
and \citet{hong2009granger} introduced a non-parametric test for non-linear Granger causality in mean based on cross-correlations.

In this paper, we focus on testing conditional mean independence with covariates
and develop a method to test (\ref{h2}) for two main reasons.
As \citet{cook2002dimension} state, regression analysis is mostly concerned with the conditional
mean of the response given the predictors, which makes testing conditional mean independence more appealing than testing conditional independence.
Further, it is very common in practice that some given covariates $Z$ have been known to affect the conditional mean of $Y$.
In this situation, we aim to determine whether $X$ has marginal effect on the conditional mean of $Y$ in the presence of $Z$,
and decide whether $X$ should be included to model the conditional mean of $Y$ along with $Z$.
In general, testing (\ref{h2}) is more useful than testing (\ref{h1}),
but requires more careful handling.

We first simplify testing (\ref{h2})
to testing conditional mean independence through a transformation.
Let $V = Y - \textrm{E}(Y|Z) \in \mathbb{R}^q$, and $U = (X, Z) \in \mathbb{R}^{p+r}$.
Then E$(V) = 0$, and E$(V|U) = \textrm{E}(Y|X,Z) - \textrm{E}(Y|Z)$.
As a result, we obtain an equivalent hypothesis test to (\ref{h2}) as
\begin{equation}\label{h3}
H_0: \textrm{E}(V|U) = \textrm{E}(V) = 0 \,\,\,\, a.s., \quad H_A: o.w.
\end{equation}
which is conditional mean independence of $V$ given $U$.
Thus, we consider the MDD with $U$ and $V$ to investigate (\ref{h3}).
However, there are two problems to solve when we apply MDD to $U$ and $V$.
First, $V$ needs to be estimated since it is unobserved.
We will replace $V$ by its estimate $\widehat{V}$ in calculating MDD.
Second, we need to confirm that the estimation error of $\widehat{V}$ is negligible,
i.e., the MDD with $\widehat{V}$ is close enough to that with $V$,
so that $\widehat{V}$ may be used for inference instead of $V$.

The rest of this paper is organized as follows.
In Section \ref{mdd}, we give a brief overview of martingale difference divergence.
In Section \ref{method}, we estimate $V$ based on the assumption that $\textrm{E}(Y|Z)$ is a linear function of $Z$,
and prove that the estimation of $V$ does not affect the asymptotic distribution of martingale difference divergence
under some regularity conditions.
We present simulation results in Section \ref{sim}, followed by a real data analysis in Section
\ref{data}\footnote{See CRAN for an accompanying R package \texttt{EDMeasure}.}. 
Finally, we summarize our work in Section \ref{con}.

The following notation is used throughout this paper.
Let $\{(X_i, Y_i, Z_i): i = 1, \dots, n\}$ be an i.i.d.\ sample
from the joint distribution $F_{X,Y,Z}$.
When $A$ is a matrix, the element of $A$ at row $k$ and column $\ell$ is denoted by $A(k, \ell)$.
When $A$ is a vector, the element of $A$ at index $k$ is denoted by $A(k)$.
The Frobenius norm of matrix $A$ $\in$ $\mathbb{R}^{p \times q}$ is denoted by $\|A\|_\textrm{F}$.
The Euclidean norm of vector $X$ $\in$ $\mathbb{R}^p$ is denoted by $|X|$.
The weighted $\mathcal{L}_2$ norm $\|\cdot\|_w$ of any complex-valued function $\eta(t), t \in \mathbb{R}^p$ is defined by $\|\eta(t)\|^2_w = \int_{\mathbb{R}^{p}} |\eta(t)|^2 w(t) \,dt$ where $|\eta(t)|^2 = \eta(t)\overline{\eta(t)}$, $\overline{\eta(t)}$ is the complex conjugate of $\eta(t)$, and $w(t)$ is any positive weight function under which the integral exists. Furthermore, $a.s.$ is an abbreviation of almost surely.

\section{MARTINGALE DIFFERENCE DIVERGENCE}\label{mdd}

\citet{shao2014martingale} proposed martingale difference divergence to capture
the conditional mean dependence (in any form) of $Y \in \mathbb{R}^q$ given $X \in \mathbb{R}^p$.

The non-negative martingale difference divergence for $X$ and $Y$, MDD$(Y | X)$ is defined by
\begin{eqnarray*}
  \textrm{MDD}^2(Y | X) = \|\textrm{E}(Y e^{i \langle s, X\rangle}) - \textrm{E}(Y)\textrm{E}(e^{i \langle s, X\rangle})\|^2_{w_p} \\
   \triangleq \int_{\mathbb{R}^{p}}{|\textrm{E}(Y e^{i \langle s, X\rangle}) - \textrm{E}(Y)\textrm{E}(e^{i \langle s, X\rangle})|^2 w_p(s) \, ds },
\end{eqnarray*}
where weight $w_p(s) = {c_p|s|^{1+p}}$, $c_p = \frac{\pi^{(1+p)/2}}{\Gamma((1+p)/2)}$, and $\Gamma$ is the gamma function.
If E$(|X|^2 + |Y|^2) < \infty$, then $\textrm{MDD}(Y | X) = 0$ if and only if $\textrm{E}(Y | X) = \textrm{E}(Y)$ holds $a.s.$.

The non-negative empirical martingale difference divergence MDD$_n(Y | X)$ is defined by
\begin{equation*}
  \textrm{MDD}_n^2(Y | X) = \frac{1}{n^2} \sum_{i,j=1}^n A_{ij} B_{ij},
\end{equation*}
where
$A_{ij} = a_{ij} - \bar{a}_{i \cdot} - \bar{a}_{\cdot j} + \bar{a}_{\cdot \cdot}$,
$\bar{a}_{i \cdot} = \frac{1}{n} \sum_{j=1}^n a_{ij}$, $\bar{a}_{\cdot j} = \frac{1}{n}\sum_{i=1}^n a_{ij}$,
$\bar{a}_{\cdot \cdot} = \frac{1}{n^2} \sum_{i,j=1}^n a_{ij}$,
$a_{ij} = |X_i - X_j|$,
and similarly for $B_{ij}$ with $b_{ij} = \frac{1}{2}|Y_i - Y_j|^2$.

The consistency and weak convergence of MDD$_n(Y | X)$ are derived as follows.
If E$(|X| + |Y|^2) < \infty$, we have
(i) $\textrm{MDD}_n(Y | X) \underset{n \rightarrow \infty}{\overset{a.s.}{\longrightarrow}} \textrm{MDD}(Y | X)$;
(ii) under $H_0: \textrm{E}(Y | X) = \textrm{E}(Y) \,\, a.s.$,
$n\textrm{MDD}_n^2(Y | X) \underset{n \rightarrow \infty}{\overset{\mathcal{D}}{\longrightarrow}} \|\zeta(s)\|^2_{w_p}$,
where $\zeta(\cdot)$ is a complex-valued zero-mean Gaussian process whose
covariance function depends on $F_{X, Y}$;
(iii) under $H_A: o.w.$,
$n\textrm{MDD}_n^2(Y | X) \underset{n \rightarrow \infty}{\overset{a.s.}{\longrightarrow}} \infty$.
Utilizing the nice properties of MDD, we next propose our test for (\ref{h3}).

\section{METHODOLOGY}\label{method}

Inspired by the linear assumption to simplify the conditional dependence structure in \citet{fan2015conditional},
we assume that the conditional expectation E$(Y | Z)$ is a linear function of $Z$,
simplifying the conditional mean dependence structure.
As a result, we can decompose $Y$ into the conditional expectation and reminder as
\begin{equation*}
Y = \textrm{E}(Y | Z) + [Y - \textrm{E}(Y | Z)] \triangleq BZ + V,
\end{equation*}
where $B \in \mathbb{R}^{q \times r}$, $V \in \mathbb{R}^q$. Then we have E$(V | Z) = 0$, and E$(V) = 0$.
Similarly, the $i$th sample counterpart is $Y_i = \textrm{E}(Y_i | Z_i) + V_i \triangleq B Z_i + V_i$, $i = 1, \dots, n$.

Suppose $\widehat{B}$ is the ordinary least squares (OLS) estimator of $B$ when regressing $Y$ on $Z$.
We will then replace $B$ with $\widehat{B}$ to estimate $\textrm{E}(Y_i | Z_i)$ as $\widehat{\textrm{E}}(Y_i | Z_i) = \widehat{B}Z_i$,
and $V_i$ as $\widehat{V}_i = Y_i - \widehat{B}Z_i = (B - \widehat{B})Z_i + V_i$.
When estimating $B$ via the OLS, $Z$ is implicitly assumed to have full column rank.
In case $Z$ is high-dimensional, i.e., $r > n$, we can estimate $B$ by the penalized least squares (PLS) similar to \citet{fan2015conditional},
including ridge \citep{hoerl1970ridge} and lasso \citep{tibshirani1996regression}.

We now construct a test for (\ref{h3}) based on $\textrm{MDD}_n^2(\widehat{V} | U)$ and its counterparts using permutation samples,
then calculate the empirical p-value following the permutation in \citet{park2015partial}.
Because the samples are independent, but with an unspecified distribution, permutation tests are a convenient tool for inference.
We will later show in Theorem 2 that the asymptotic distribution of $n\textrm{MDD}_n^2(\widehat{V} | U)$ depends on an unknown underlying distribution,
which justifies the use of permutation tests.
To measure the conditional mean dependence of $V$ given $U$,
we first compute the test statistic $\textrm{MDD}_n^2(\widehat{V} | U)$ from
the sample $\{(\widehat{V}_i, U_i): i = 1, \dots, n\}$,
where $U_i = (X_i, Z_i)$. That is, $\textrm{MDD}_n^2(\widehat{V} | U)$ depends on the i.i.d.\ sample $\{(X_i, Y_i, Z_i): i = 1, \dots, n\}$.
Next we draw $B$ permutation samples of size $n$ as $\{(X_i^\ast, Y_i, Z_i): i = 1, \cdots, n\}$,
where only the sample of $X$ is permuted in order to approximate the sampling distribution. For each permutation sample, we calculate the
test statistic $\textrm{MDD}_{n,b}^2(\widehat{V} | U)$, $b = 1, \cdots, B$.
Then the empirical p-value is given by
\begin{equation*}
\widehat{p} = \frac{\sum_{b=1}^B \mathbf{1}\left\{\textrm{MDD}_{n,b}^2(\widehat{V} | U) \geq \textrm{MDD}_n^2(\widehat{V} | U)\right\}}{B}.
\end{equation*}

When $H_0$ is false, $\textrm{MDD}_n^2(\widehat{V} | U)$ tends to be large while $\textrm{MDD}_{n,b}^2(\widehat{V} | U)$ tends to be small.
As a result, the empirical p-value is expected to be very small, leading to a rejection of $H_0$.
We name the proposed test linear martingale difference divergence (LinMDD). To justify our LinMDD test, it remains to validate that
$\textrm{MDD}_n^2(\widehat{V}|U)$ is close enough to $\textrm{MDD}_n^2(V|U)$,
i.e., the estimation error in $\widehat{V}$ is negligible for the sampling distribution
of the test statistic, focusing on the asymptotic case.
To begin with, we introduce some regularity conditions to derive the asymptotic distribution of $\textrm{MDD}_n^2(\widehat{V}|U)$.

\begin{condition}
There exist constants $0 < c_1, c_2, c_3 < \infty$, such that
$\textrm{E}(|U_i - U_j|^2) = c_1$, $i \neq j$;
$\textrm{E}(|U_i - U_j||U_i - U_k|) = c_2$, $i \neq j \neq k$;
$\textrm{E}(|U_i - U_j||U_k - U_\ell|) = c_3$, $i \neq j \neq k \neq \ell$.
\end{condition}
\begin{condition}
There exists constant $0 < c_4 < \infty$, such that
$\textrm{E}[(Z_i(t) - Z_j(t))^2 (Z_i(s) - Z_j(s))^2] \leq c_4$, $i \neq j$, $\forall t, s$.
\end{condition}
\begin{condition}
There exists constant $0 < c_5 < \infty$, such that
$\textrm{E}[(Z_i(t) - Z_j(t))^2 (V_i(s) - V_j(s))^2] \leq c_5$, $i \neq j$, $\forall t, s$.
\end{condition}
\begin{condition}
$\|\widehat{B} - B\|_F = O_p(n^{-1/2})$.
\end{condition}

\begin{remark}
Condition 4 can be derived from the bounded density of $|V_i - V_j|$ and non-heavy tails of $Z_i(t)$ and $V_i(t)$
according to \citet{fan2015conditional} and \citet{fan2011high}.
\end{remark}

Through a similar derivation to Theorem 2 of \citet{fan2015conditional},
we justify the choice of using $\textrm{MDD}_n^2(\widehat{V}|U)$ in place of $\textrm{MDD}_n^2(V|U)$ by the following lemma and theorems.
Lemma \ref{conv} shows that the difference between $\textrm{MDD}_n^2(\widehat{V} | U)$ and $\textrm{MDD}_n^2(V | U)$ is negligible as the sample size increases.
The proof of Lemma \ref{conv} can be found in Appendix \ref{pf_lemma1}.

\begin{lemma}\label{conv}
If $Y = BZ + V$ and Conditions 1-4 hold, we have
\begin{equation*}
\textrm{MDD}_n^2(\widehat{V} | U) - \textrm{MDD}_n^2(V | U) = O_p(n^{-3/2}).
\end{equation*}
\end{lemma}

Consequently, the consistency and weak convergence of MDD$_n(\widehat{V} | U)$ follow from Lemma \ref{conv}
and are summarized in Theorem \ref{thm1} and \ref{thm2} below.
\begin{theorem}[Consistency]\label{thm1}
If $Y = BZ + V$ and Conditions 1-4 hold, we have
\begin{equation*}
\textrm{MDD}_n(\widehat{V} | U) \underset{n \rightarrow \infty}{\overset{\mathcal{P}}{\longrightarrow}} \textrm{MDD}(V | U).
\end{equation*}
\end{theorem}

\begin{theorem}[Weak convergence]\label{thm2}
If $Y = BZ + V$ and Conditions 1-4 hold, under $H_0$, we have
\begin{equation*}
n\textrm{MDD}_n^2(\widehat{V} | U) \underset{n \rightarrow \infty}{\overset{\mathcal{D}}{\longrightarrow}} \|\zeta(s)\|^2_{w_p},
\end{equation*}
where $\zeta(\cdot)$ denotes the complex-valued Gaussian random process corresponding to
the asymptotic distribution of $n\textrm{MDD}_n^2(V | U)$.
Under $H_A$, we have
\begin{equation*}
n\textrm{MDD}_n^2(\widehat{V} | U) \underset{n \rightarrow \infty}{\overset{\mathcal{P}}{\longrightarrow}} \infty.
\end{equation*}
\end{theorem}

According to Theorem \ref{thm1},
$\textrm{MDD}_n(\widehat{V} | U)$ converges to the same population statistic $\textrm{MDD}(V | U)$ as $\textrm{MDD}_n(V | U)$,
and thus it can serve to measure the conditional mean dependence of $V$ given $U$.
In addition, $n\textrm{MDD}_n^2(\widehat{V} | U)$ and $n\textrm{MDD}_n^2(V | U)$
have the same asymptotic distribution stated in Theorem \ref{thm2},
which establishes the effectiveness of LinMDD test,
as we approximate the limiting distribution of $n\textrm{MDD}_n^2(V | U)$ using $n\textrm{MDD}_n^2(\widehat{V} | U)$ in LinMDD test.
In Section \ref{sim} and Section \ref{data}, we will present the finite-sample performance of our LinMDD test through simulations and a real data example, respectively.

\section{SIMULATION STUDIES}\label{sim}

To evaluate the performance of our LinMDD test,
we adopt the simulation setup in \citet{lavergne2000nonparametric},
and compare our test to the pMDD test \citep{park2015partial},
pdCov test \citep{szekely2014partial},
and pcov test \citep{lan2014testing} as benchmarks.
All tests are implemented as permutation tests with permutation size $B = 500$,
in which we only permute the sample of $X$ to approximate the distribution of the test statistic.

We generate data from the underlying model
\begin{equation*}
Y = - Z + b \cdot Z^3 + f(X) + \epsilon,
\end{equation*}
where $Z \sim N(0, 1)$, $X \sim \mathcal{N}(0, 1)$, $\epsilon \sim \mathcal{N}(0, 4)$, and $Z, X, \epsilon$ are independent.
We test the null hypothesis $H_0: \textrm{E}(Y | X, Z) = \textrm{E}(Y | Z) \,\, a.s.$ with significance level $\alpha \in \{0.05, 0.1\}$,
and examine the empirical size and power of each test.
We run 1000 replications with sample size $n \in \{20, 30, 50, 70, 100\}$ for each specific model.

\begin{experiment}[Linear $Z$, linear $X$]\label{exp1}
$b = 0$, $f(X) = cX$ where $c \in \{0, \frac{2}{3}, 1, \frac{3}{2}\}$.
\end{experiment}

\begin{experiment}[Linear $Z$, non-linear $X$]\label{exp2}
$b = 0$, $f(X) = \sin(c \pi X)$ where $c \in \{\frac{1}{4}, \frac{1}{3}, \frac{1}{2}\}$.
We omit $c = 0$ as it is exactly the same as $c = 0$ in Model \ref{exp1}.
\end{experiment}

\begin{figure}[ht]
\begin{center}
\centerline{\includegraphics[width=\columnwidth]{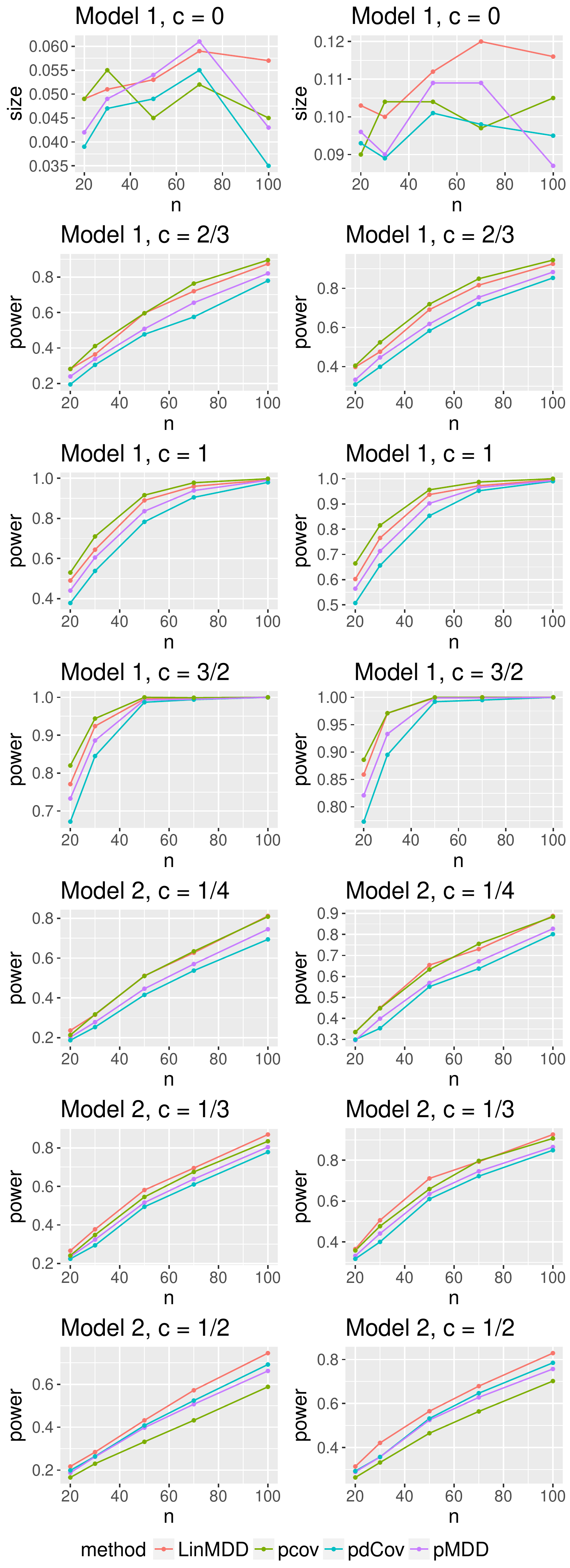}}
\caption{Empirical size and power of 1000 replications with $B = 500$ for Model \ref{exp1} \& \ref{exp2}.}\label{fig1}
\end{center}
\vskip -0.2in
\end{figure}

From Figure \ref{fig1}, the empirical size of all tests is around 0.05 (0.1).
The empirical power of all tests increases as $n$ increases.
For the linear $X$ case, the empirical power of all tests is higher when $c$ is larger,
since the signal-to-noise ratio increases.
Moreover,
the empirical power of the LinMDD and pcov tests is consistently higher than that of the other tests,
because the linear assumption is valid, and only LinMDD and pcov tests are designed for linear $Z$.
For the non-linear $X$ case,
the LinMDD test still outperforms the other tests,
while the performance of the pcov test degrades as $c$ increases,
because the LinMDD test is designed for non-linear $X$ while pcov test is suitable only for linear $X$.

\begin{experiment}[Nonlinear $Z$, linear $X$]\label{exp3}
$b = 1$, $f(X) = cX$ where $c \in \{0, \frac{2}{3}, 1, \frac{3}{2}\}$.
\end{experiment}

\begin{experiment}[Nonlinear $Z$, non-linear $X$]\label{exp4}
$b = 1$, $f(X) = \sin(c \pi X)$ where $c \in \{\frac{1}{4}, \frac{1}{3}, \frac{1}{2}\}$.
We omit $c = 0$ as it is exactly the same as $c = 0$ in Model \ref{exp3}.
\end{experiment}

\begin{figure}[ht]
\begin{center}
\centerline{\includegraphics[width=\columnwidth]{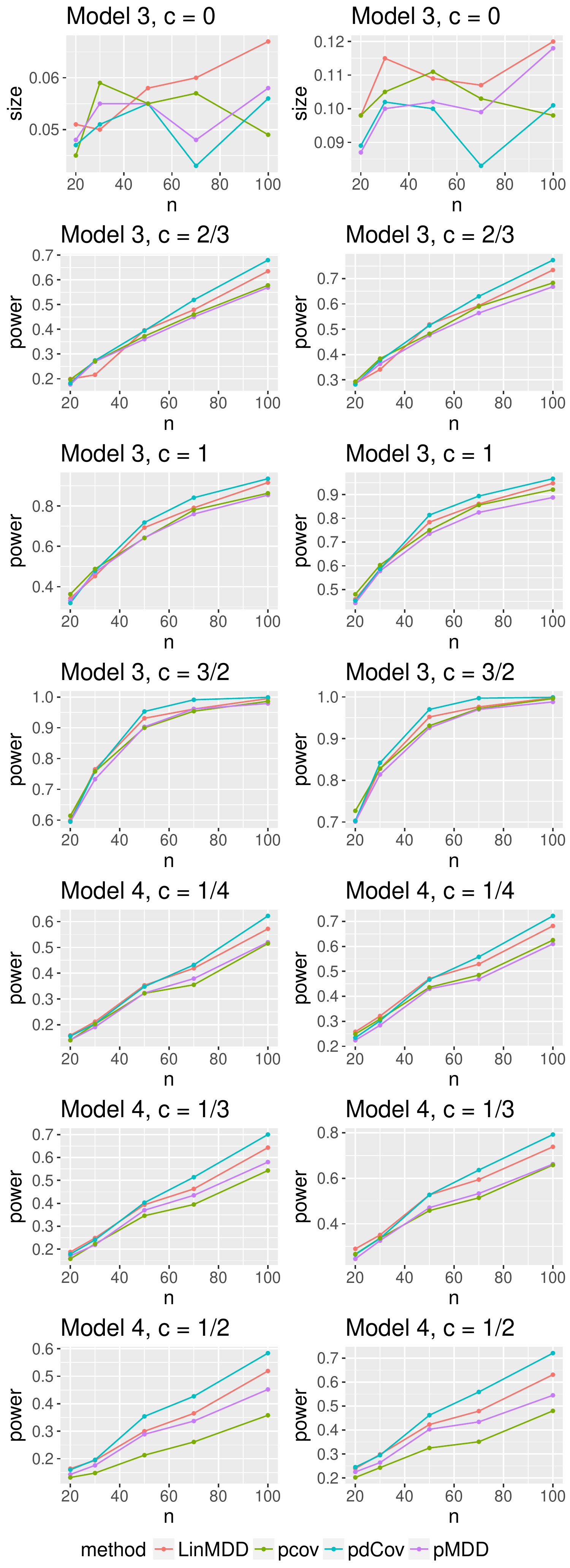}}
\caption{Empirical size and power of 1000 replications with $B = 500$ for Model \ref{exp3} \& \ref{exp4}.}\label{fig2}
\end{center}
\vskip -0.2in
\end{figure}

From Figure \ref{fig2}, 
the empirical size of all tests is around 0.05 (0.1).
For the linear $X$ case,
the empirical power of the LinMDD and pcov tests is competitive with but not always higher than that of the other tests.
The reason is that the linear dependence of $Y$ on $Z$ is violated
while the other tests do not rely it.
For the non-linear $X$ case, we similarly find that the performance of the pcov test degrades as $c$ increases.
The simulation results show that our LinMDD test achieves competitive and often better performance than
the others in these situations. Next, we apply the proposed LinMDD test on a real dataset.

\section{FINANCIAL DATA APPLICATION}\label{data}

In finance, the capital asset pricing model (CAPM) was proposed by
\citet{sharpe1964capital}, \citet{lintner1965valuation}, and \citet{mossin1966equilibrium}
to describe the stock returns through the market risk as
\begin{equation*}
r_t = \alpha + \beta_1 m_t,
\end{equation*}
where $r_t$ is the excess stock return (in excess the risk-free return),
and $m_t$ is the excess market return at time $t$.
\citet{fama1993common} added size and value factors to the CAPM, and proposed the Fama$-$French three-factor model as
\begin{equation*}
r_t = \alpha + \beta_1 m_t + \beta_2 \, \textrm{SMB}_t + \beta_3 \, \textrm{HML}_t,
\end{equation*}
where SMB (small minus big) and HML (high minus low) account for
stocks with small/big market capitalization and high/low book-to-market ratio, respectively.
\citet{fama2015five} further added profitability and investment factors to the three-factor model, and extended it to the Fama$-$French five-factor model as
\begin{eqnarray*}
r_t = \alpha + \beta_1 m_t + \beta_2 \, \textrm{SMB}_t + \beta_3 \, \textrm{HML}_t \\
+ \,\beta_4 \, \textrm{RMW}_t + \beta_5 \, \textrm{CMA}_t,
\end{eqnarray*}
where RMW (robust minus weak) and CMA (conservative minus aggressive) further account for
stocks with robust/weak operating profitability and conservative/aggressive investment, respectively.

\begin{figure*}[ht]
\begin{center}
\centerline{\includegraphics[width=1.6\columnwidth]{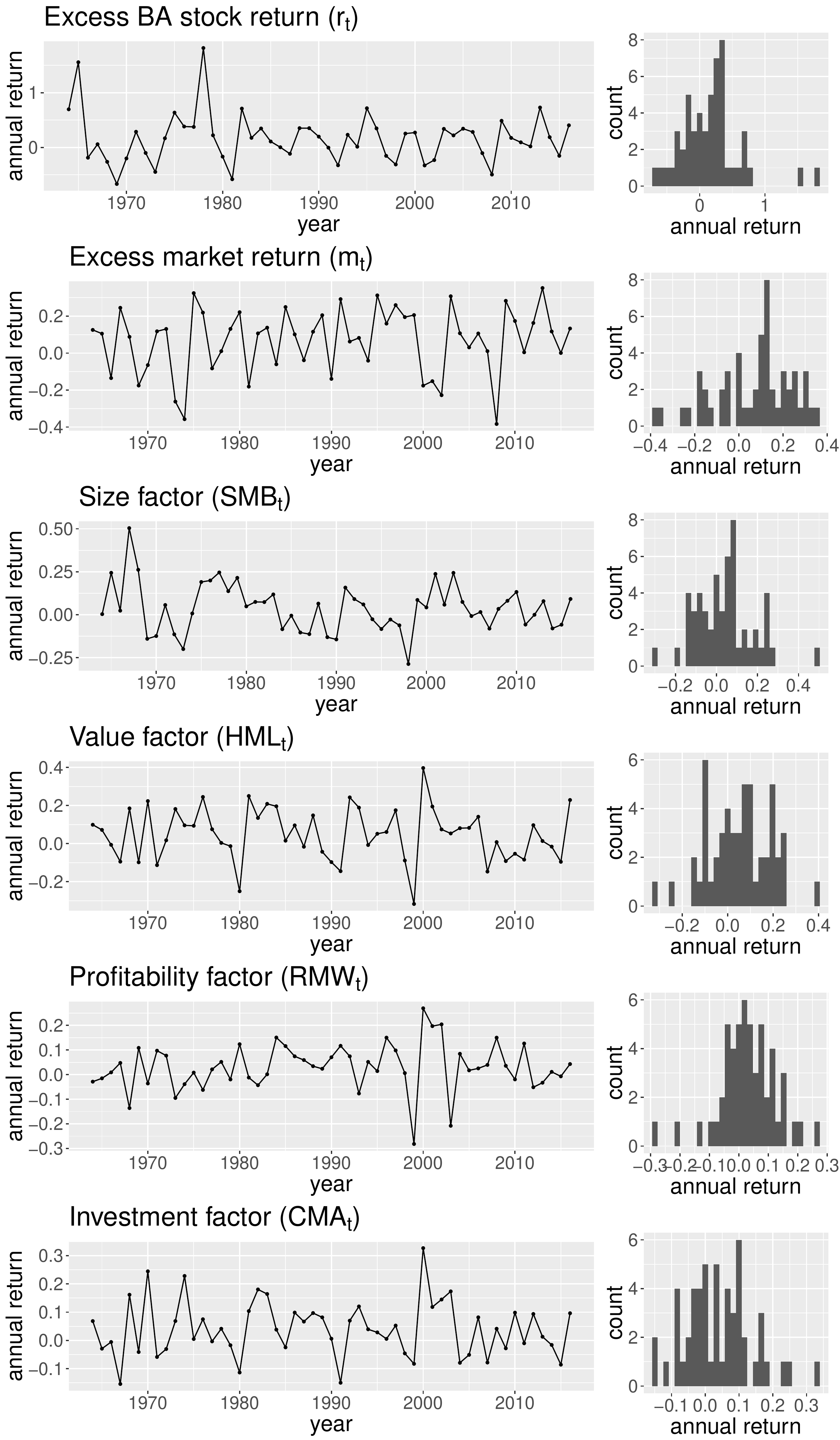}}
\caption{Time series and histograms of excess BA stock returns ($r_t$),
excess market returns ($m_t$), size factors ($\textrm{SMB}_t$),
value factors ($\textrm{HML}_t$), profitability factors ($\textrm{RMW}_t$),
and investment factors ($\textrm{CMA}_t$) between 1964 and 2016.}\label{fig3}
\end{center}
\vskip -0.2in
\end{figure*}

We collect the annual risk-free returns and
Fama$-$French five factors\footnote{Download data at http://mba.tuck.dartmouth.edu/pages/fa\\culty/ken.french/data\_library.html.},
and the annual returns of Boeing (BA) stock\footnote{Download data using \texttt{get.hist.quote} in the R package \texttt{tseries} \citep{adrian2017tseries}.}
in the past 53 years between 1964 and 2016.
The time series and histograms of excessive BA stock returns and Fama$-$French five factors
are depicted in Figure \ref{fig3}.

\subsection{CAPM VS. FAMA$-$FRENCH THREE-FACTOR MODEL}\label{data1}

First, we are curious whether the size and value factors should be added to the CAPM,
i.e., whether SMB and HML in the Fama$-$French three-factor model contribute to the expectation of excess stock returns given the market risk.
Thus, we test $H_0: \textrm{E}(Y | X, Z) = \textrm{E}(Y | Z) \,\, a.s.$,
where $X_t = (\textrm{SMB}_t, \textrm{HML}_t)$, $Y_t = r_t$, and $Z_t = (1, m_t)$.

We apply our LinMDD test to the data with $n = 53$ and $B = 500$. Our p-value is 0.072,
while the p-values are 0.012 (pMDD) and 0.092 (pdCov) using competing tests.
As a result, we reject $H_0$ with significance level $\alpha = 0.1$, and conclude that
SMB and HML help determine the excess returns of BA stock in the presence of the market risk.
Our results align with the research in finance that
the Fama$-$French three-factor model remarkably outperforms the CAPM in explaining excess stock returns.

\subsection{FAMA$-$FRENCH THREE-FACTOR VS. FIVE-FACTOR MODEL}\label{data2}

Similarly, we are interested in whether the profitability and investment factors should be further added to the Fama$-$French three-factor model,
i.e., whether RMW and CMA in the Fama$-$French five-factor model contribute to the description of excess stock returns given the other three factors.
Hence, we test $H_0: \textrm{E}(Y | X, Z) = \textrm{E}(Y | Z) \,\, a.s.$,
in which $X_t = (\textrm{RMW}_t, \textrm{CMA}_t)$, $Y_t = r_t$, and $Z_t = (1, m_t, \textrm{SMB}_t, \textrm{HML}_t)$.

We apply our LinMDD test to the data with $n = 53$ and $B = 500$, and its p-value is 0.360,
while the p-values are 0.358 (pMDD) and 0.878 (pdCov) using competing tests.
As a result, we fail to reject $H_0$ with significance level $\alpha = 0.1$, and conclude that
RMW and CMA are unable to help determine the excess returns of BA stock in the presence of the other three factors.
Our results align with the research in finance that
the Fama$-$French five-factor model has yet to be proven as a significant improvement over the three-factor model in describing excess stock returns.

\subsection{FAMA$-$FRENCH FOUR-FACTOR MODEL VS. FIVE-FACTOR MODEL}\label{data3}

\citet{fama2015five} showed that the value factor HML becomes redundant
when profitability and investment factors are added to the Fama$-$French three-factor model,
because HML is fully captured by its exposures to the other four factors, especially RMW and CMA.
To validate this argument, we test $H_0: \textrm{E}(Y | X, Z) = \textrm{E}(Y | Z) \,\, a.s.$,
where $X_t = \textrm{HML}_t$, $Y_t = r_t$, and $Z_t = (1, m_t, \textrm{SMB}_t, \textrm{RMW}_t, \textrm{CMA}_t)$.

We apply our LinMDD test to the data with $n = 53$ and $B = 500$. Our p-value is 0.218,
while the p-values are 0.438 (pMDD) and 0.540 (pdCov) using competing tests.
As a result, we fail to reject $H_0$ with significance level $\alpha = 0.1$, and conclude that
HML cannot help explain the excess returns of BA stock in the presence of the other four factors.
Our results demonstrate that HML is redundant for describing excess stock returns in the Fama$-$French five-factor model.

\section{CONCLUSION}\label{con}

In this paper, we propose a new test, LinMDD,
for the null hypothesis $H_0: \textrm{E}(Y | X, Z) = \textrm{E}(Y | Z) \,\, a.s.$
by investigating an equivalent one $H_0: \textrm{E}(V|U) = \textrm{E}(V) = 0 \,\, a.s.$,
derived from a transformation involving the conditional expectation.
When applying martingale difference divergence \citep{shao2014martingale} to
test $H_0: \textrm{E}(V|U) = \textrm{E}(V) = 0 \,\, a.s.$, we make two major contributions.

(1) Since $V$ is unobservable, we estimate $V$ based on the assumption that
$\textrm{E}(Y | Z)$ is a linear function of $Z$,
simplifying the conditional mean dependence structure.

(2) We prove that the estimation error in $\widehat{V}$ is negligible for the asymptotic distribution
of the test statistic. Thus, we can replace $V$ with $\widehat{V}$ in the test statistic for inference in large samples.

We implement the LinMDD test as a permutation test following \citet{park2015partial},
and compare it with existing tests in various simulation studies.
The LinMDD test consistently outperforms existing tests when its linear assumption is valid, and
it achieves competitive results with existing tests even when its linear assumption is violated.

To illustrate the practical value of the LinMDD test,
we compare the CAPM, the Fama$-$French three-factor and five-factor models by applying LinMDD test to the financial data.
We find that the Fama$-$French three-factor outperforms the CAPM,
while the Fama$-$French five-factor is not a significant improvement over the three-factor model
when explaining the excess annual returns of a major stock.
Moreover, we validate the statement that the value factor is redundant in the Fama$-$French five-factor model \citep{fama2015five}
using the LinMDD test.

The relaxation of the linear assumption is an important topic for future research.
Our method will become more general if the linear assumption of conditional mean dependence
can be generalized to a non-linear one, using non-parametric regression (local regression, splines)
instead of linear regression in the estimation of conditional mean.
In addition, the high-dimensional setting regarding $Z$ where $r>n$ is an interesting direction to consider as well.






%


%
%


\appendix
\section*{Appendix}
\addcontentsline{toc}{section}{Appendices}
\renewcommand{\thesubsection}{\Alph{subsection}}

\subsection{Proof of Lemma \ref{conv}}\label{pf_lemma1}

\begin{proof}
We define $T$
\begin{eqnarray*}
\begin{aligned}
&= n\textrm{MDD}_n^2(\widehat{V} | U) - n\textrm{MDD}_n^2(V | U) \\
&= \frac{1}{2n} \sum_{i,j} [(F_{ij}-\frac{1}{n}\sum_{k}F_{kj}-\frac{1}{n}\sum_{k}F_{ik}+\frac{1}{n^2}\sum_{k,\ell}F_{k\ell}) \\
& \times \, (E_{ij}-\frac{1}{n}\sum_{k}E_{kj}-
\frac{1}{n}\sum_{k}E_{ik}+\frac{1}{n^2}\sum_{k,\ell}E_{k\ell})],
\end{aligned}
\end{eqnarray*}
where $F_{ij} = |\widehat{V}_i - \widehat{V}_j|^2 - |V_i - V_j|^2$, $E_{ij} = |U_i - U_j|$.

We apply Taylor expansion to $|\widehat{V}_t - \widehat{V}_s|^2$ at $V_t - V_s$ in terms of $f(x) = x^Tx$, $f'(x) = 2x^T$,
then there exists $\lambda \in (0, 1)$, such that $F_{ij}$
\begin{eqnarray*}
\begin{aligned}
 &= 2[\lambda (\widehat{V}_i - \widehat{V}_j) + (1-\lambda)(V_i - V_j) ]^T (\widehat{V}_i - \widehat{V}_j - V_i + V_j) \\
 &= 2[ \lambda (Z_i - Z_j)^T(B - \widehat{B})^T(B - \widehat{B})(Z_i - Z_j) \\
 &+ \, (V_i - V_j)^T(B - \widehat{B})(Z_i - Z_j) ].
\end{aligned}
\end{eqnarray*}

Thus, we have $T = T_1 + T_2$, where $T_1$
\begin{eqnarray*}
\begin{aligned}
 &= \frac{\lambda}{n} \sum_{i,j} [(G_{ij}-\frac{1}{n}\sum_{k}G_{kj}-\frac{1}{n}\sum_{k}G_{ik}+\frac{1}{n^2}\sum_{k,\ell}G_{k\ell})\\
 &\times \, (E_{ij}-\frac{1}{n}\sum_{k}E_{kj}-\frac{1}{n}\sum_{k}E_{ik}+\frac{1}{n^2}\sum_{k,\ell}E_{k\ell})],
\end{aligned}
\end{eqnarray*}
\begin{equation*}
G_{ij} = (Z_i - Z_j)^T(B - \widehat{B})^T(B - \widehat{B})(Z_i - Z_j),
\end{equation*}
and $T_2$
\begin{eqnarray*}
\begin{aligned}
&= \frac{1}{n} \sum_{i,j} [(H_{ij}-\frac{1}{n}\sum_{k}H_{kj}-\frac{1}{n}\sum_{k}H_{ik} \frac{1}{n^2}\sum_{k,\ell}H_{k\ell}) \\
 &\times \, (E_{ij} - \frac{1}{n}\sum_{k}E_{kj}-\frac{1}{n}\sum_{k}E_{ik}+\frac{1}{n^2}\sum_{k,\ell}E_{k\ell})],
\end{aligned}
\end{eqnarray*}
\begin{equation*}
H_{ij} = (V_i - V_j)^T(B - \widehat{B})(Z_i - Z_j).
\end{equation*}

First, we will show (i) $T_1 = O_p(n^{-1})$.

After a simple calculation, we have
\begin{eqnarray*}
\begin{aligned}
&\frac{1}{n} \sum_{i,j} (G_{ij}-\frac{1}{n}\sum_{k}G_{kj}-\frac{1}{n}\sum_{k}G_{ik}+\frac{1}{n^2}\sum_{k,\ell}G_{k\ell})E_{ij} \\
&= tr[\frac{1}{n} \sum_{i,j}|U_i - U_j| (G_{ij}-\frac{1}{n}\sum_{k}G_{kj}-\frac{1}{n}\sum_{k}G_{ik} \\
& + \, \frac{1}{n^2}\sum_{k,\ell}G_{k\ell})] \\
&= tr[(B - \widehat{B})^T(B - \widehat{B}) M],
\end{aligned}
\end{eqnarray*}
where $M = \frac{1}{n} \sum_{i,j}|U_i - U_j| S_{ij}$, and
\begin{equation*}
S_{ij} = R_{ij}-\frac{1}{n}\sum_{k}R_{kj}-\frac{1}{n}\sum_{k}R_{ik}+\frac{1}{n^2}\sum_{k,\ell}R_{k\ell},
\end{equation*}
$R_{ij} = (Z_i - Z_j)(Z_i - Z_j)^T$, $R_{ij} = R_{ji}$, $S_{ij} = S_{ji}$, then
\begin{eqnarray*}
\begin{aligned}
&\textrm{E}[(M(t,s))^2 ] \\
&= \textrm{E}[\frac{1}{n^2} (\sum_{i,j}|U_i - U_j| S_{ij}(t,s))^2 ] \\
&= \textrm{E}\{ \textrm{E} [\frac{1}{n^2} (\sum_{i,j}|U_i - U_j| S_{ij}(t,s))^2 | U_i, \forall i ] \} \\
&= \textrm{E}[ \frac{2c_1}{n^2} \sum_{i\neq j} (S_{ij}(t,s))^2 \\
& + \, \frac{2c_2}{n^2} \sum_{i\neq j \neq k} (S_{ij}(t,s)S_{i k}(t,s) + S_{ij}(t,s)S_{kj}(t,s)) \\
& + \, \frac{c_3}{n^2} \sum_{i\neq j \neq k \neq \ell} S_{ij}(t,s)S_{k \ell}(t,s)],
\end{aligned}
\end{eqnarray*}
where $c_1 = \textrm{E}(|U_i - U_j|^2)$, $i \neq j$; $c_2 = \textrm{E}(|U_i - U_j||U_i - U_k|)$, $i \neq j \neq k$;
$c_3 = \textrm{E}(|U_i - U_j||U_k - U_\ell|)$, $i \neq j \neq k \neq \ell$.

Considering that $\textrm{E}[(R_{ij}(t,s))^2] = \textrm{E}[(Z_i - Z_j)_t^2 (Z_i - Z_j)_s^2] \leq c_4$, $i \neq j$, $\forall t, s$, we have
$\textrm{E}[(R_{ij}(t,s))^2] = O(1)$,
which implies $\textrm{E}[(S_{ij}(t,s))^2] = O(1)$, and thus $\textrm{E}[\frac{1}{n^2} \sum_{i\neq j} (S_{ij}(t,s))^2] = O(1)$.

After a simple calculation, we have $\sum_i S_{ij}(t,s) = 0$, $\sum_j S_{ij}(t,s) = 0$, $\sum_i\sum_j S_{ij}(t,s) = 0$, and
\begin{eqnarray*}
\begin{aligned}
& \sum_{i\neq j \neq k} S_{ij}(t,s)S_{i k}(t,s) \\
&= \sum_i (S_{ii}(t,s))^2  - \sum_{i \neq j} (S_{ij}(t,s))^2, \\
& \sum_{i\neq j \neq k} S_{ii}(t,s)S_{jk}(t,s) \\
&= \sum_i (S_{ii}(t,s))^2 - \sum_{i \neq j} S_{ii}(t,s) S_{j j}(t,s), \\
& \sum_{i\neq j \neq k \neq \ell} S_{ij}(t,s)S_{k \ell}(t,s) \\
&= - \, 2\sum_{i \neq j \neq k} [S_{ii}(t,s)S_{jk}(t,s) + S_{ij}(t,s)S_{ik}(t,s) \\
& + \, S_{ij}(t,s)S_{kj}(t,s)] \\
& - \, \sum_{i \neq j} [4S_{ii}(t,s)S_{ij}(t,s) + S_{ii}(t,s)S_{jj}(t,s) \\
& + \, 2(S_{ij}(t,s))^2] - \sum_{i} (S_{ii}(t,s))^2,
\end{aligned}
\end{eqnarray*}
we have
\begin{eqnarray*}
\begin{aligned}
\textrm{E}[\frac{1}{n^2} \sum_{i\neq j \neq k} S_{ij}(t,s)S_{i k}(t,s)] &= O(1), \\
\textrm{E}[\frac{1}{n^2} \sum_{i\neq j \neq k} S_{ij}(t,s)S_{kj}(t,s)] &= O(1), \\
\textrm{E}[\frac{1}{n^2} \sum_{i\neq j \neq k} S_{ii}(t,s)S_{jk}(t,s)] &= O(1), \\
\textrm{E}[\frac{1}{n^2} \sum_{i\neq j \neq k \neq \ell} S_{ij}(t,s)S_{k \ell}(t,s)] &= O(1).
\end{aligned}
\end{eqnarray*}
Therefore, $\textrm{E}[(M(t,s))^2] = O(1)$.

Applying Chebyshev's inequality to $M(t,s)$, we have
\begin{equation*}
P(|M(t,s) - \mu| \geq k\sigma) \leq 1/k^2,
\end{equation*}
where $\mu = \textrm{E}[M(t,s)]$, $\sigma^2 = \textrm{Var}[M(t,s)]$.
As a result, $M(t,s) = O_p(1)$.

Given that $\|\widehat{B} - B\|_F = O_p(n^{-1/2})$, we have
\begin{eqnarray*}
\begin{aligned}
& \frac{1}{n} \sum_{i,j} (G_{ij}-\frac{1}{n}\sum_{k}G_{kj}-\frac{1}{n}\sum_{k}G_{ik}+\frac{1}{n^2}\sum_{k,\ell}G_{k\ell})E_{ij} \\
&= tr[(B - \widehat{B})^T(B - \widehat{B}) M] \\
&= pq^2 O_p(n^{-1}) O_p(1) \\
&= O_p(n^{-1}).
\end{aligned}
\end{eqnarray*}

Similarly, we have
\begin{eqnarray*}
\begin{aligned}
&\frac{1}{n} \sum_{i,j} (G_{ij}-\frac{1}{n}\sum_{k}G_{kj}-\frac{1}{n}\sum_{k}G_{ik}+\frac{1}{n^2}\sum_{k,\ell}G_{k\ell})E_{kj}, \\
&\frac{1}{n} \sum_{i,j} (G_{ij}-\frac{1}{n}\sum_{k}G_{kj}-\frac{1}{n}\sum_{k}G_{ik}+\frac{1}{n^2}\sum_{k,\ell}G_{k\ell})E_{ik}, \\
&\frac{1}{n} \sum_{i,j} (G_{ij}-\frac{1}{n}\sum_{k}G_{kj}-\frac{1}{n}\sum_{k}G_{ik}+\frac{1}{n^2}\sum_{k,\ell}G_{k\ell})E_{k \ell}
\end{aligned}
\end{eqnarray*}
are all $O_p(n^{-1})$. Therefore, $T_1 = O_p(n^{-1})$.

Analogous to (i), we can show (ii) $T_2 = O_p(n^{-1/2})$.
The only differences are
\begin{eqnarray*}
\begin{aligned}
& \frac{1}{n} \sum_{i,j} (H_{ij}-\frac{1}{n}\sum_{k}H_{kj}-\frac{1}{n}\sum_{k}H_{ik}+\frac{1}{n^2}\sum_{k,\ell}H_{k\ell})E_{ij} \\
&= tr[(B - \widehat{B}) M],
\end{aligned}
\end{eqnarray*}
where $M$ is defined similarly with
$R_{ij} = (Z_i - Z_j)(V_i - V_j)^T$, and
$\textrm{E}[(R_{ij}(t,s))^2] = \textrm{E}[(Z_i - Z_j)_t^2 (V_i - V_j)_s^2] \leq c_5$, $i \neq j$, $\forall t, s$, and
\begin{eqnarray*}
\begin{aligned}
& \frac{1}{n} \sum_{i,j} (H_{ij}-\frac{1}{n}\sum_{k}H_{kj}-\frac{1}{n}\sum_{k}H_{ik}+\frac{1}{n^2}\sum_{k,\ell}H_{k\ell})E_{ij} \\
&= tr[(B - \widehat{B}) M] \\
&= pq O_p(n^{-1/2}) O_p(1) \\
&= O_p(n^{-1/2}),
\end{aligned}
\end{eqnarray*}
and therefore $T_2 = O_p(n^{-1/2})$.

As a conclusion, $T = T_1 + T_2 = O_p(n^{-1/2})$.
\end{proof}

\clearpage

\renewcommand{\refname}{\normalsize{References}}
\bibliographystyle{abbrvnat}
\bibliography{Paper}

\end{document}